\documentclass[amsppt,11pt]{amsart}

\usepackage{amsmath}
\usepackage{amsthm}
\usepackage{amssymb}
\usepackage{verbatim}

\allowdisplaybreaks[4]

\def\Z{\mathbb {Z}}

\def\E{\mathbb{E}}

\def\P{\mathbb{P}}

\def\1{{\mathbf 1}}
\def\0{{\mathbf 0}}

\newtheorem{theorem}{Theorem}[section]
\newtheorem{lemma}{Lemma}[section]

\def\var{\mathrm {Var}}

\def\Cov{\mathrm {Cov}}

\newcommand{\dir}[2]{\mathcal{E} (#1\mid #2)}
\newcommand{\dirt}[2]{\mathcal{E}_t (#1 \mid #2)}

\usepackage{verbatim, epsfig}
\title{Mixing times for random walks on finite lamplighter groups}
\author[Yuval Peres\,\, David Revelle]
{Yuval Peres$^*$, David Revelle$^\dagger$}
\date{March 29, 2004.
\newline\indent
$^*$Research supported by  NSF Grants
DMS-0104073 and DMS-0244479.
\newline\indent
$^\dagger$Research partially supported by an NSF postdoctoral
fellowship}
\def\lamp{G_n^\diamondsuit}
\def\lampp{G^\diamondsuit}
\def\lampzd{(\Z_n^d)^\diamondsuit}
\def\lampztwo{(\Z_n^2)^\diamondsuit}
\def\Trel{T_{\mathrm {rel}}}
\def\Cn{\mathcal{C}_n}
\def\C{\mathcal{C}}

\def\Tret{T_x^+}
\newcommand{\T}{{\mathcal T}_{\mathrm v}}
\newcommand{\Tg}[1]{\T(#1)}
\newcommand{\Tve}[1]{\T(\epsilon, #1)}
\newcommand{\Tvef}[1]{\T(\epsilon/2, #1)}
\newcommand{\Teg}[1]{\T(1/(2e),#1)}

\newcommand{\tauve}[1]{\tau(\epsilon, #1)}
\newcommand{\tauvef}[1]{\tau(\epsilon/2, #1)}
\newcommand{\tauveg}[2]{\tau(#1,#2)}
\newcommand{\taug}[1]{\tau(1/2, #1)}

\begin{document}

\begin{abstract}
  Given a  finite graph $G$, a vertex of the lamplighter
  graph $\lampp$ consists of a zero-one labeling of the vertices of
  $G$, and a marked vertex of $G$.  If $G$ is a Cayley graph, then
  $\lampp$ is the wreath product $\Z_2 \wr G$. 
  For transitive $G$ we show that, up to constants, the relaxation
  time for simple random walk in $\lampp$ is the maximal hitting time
  for simple random walk in $G$, while the mixing time in total
  variation on $\lampp$ is the expected cover time on $G$.  The mixing
  time in the uniform metric on $\lampp$ admits a sharp threshold, and
  equals $|G|$ multiplied by the relaxation time on $G$, up to a
  factor of $\log |G|$.  For $\Z_2 \wr \Z_n^2$, the lamplighter group
  over the discrete two dimensional torus, the relaxation time is of
  order $n^2 \log n$, the total variation mixing time is of order $n^2
  \log^2 n$, and the uniform mixing time is of order $n^4$.  For $\Z_2
  \wr \Z_n^d$ when $d\geq 3$, the relaxation time is of order $n^d$,
  the total variation mixing time is of order $n^d \log n$, and the
  uniform mixing time is of order $n^{d+2}$. These are the first
  examples we know of finite transitive graphs with uniformly bounded
  degrees where these three mixing time parameters are of different
  orders of magnitude.
    \end{abstract}

\maketitle

\section{Introduction}
Given a finite graph $G=(V_G,E_G)$, the wreath product $\lampp=\Z_2
\wr G$ is the graph whose vertices are ordered pairs $(f,x)$, where $x
\in V_G$ and $f\in \{0,1\}^{V_G}$.  There is an edge between $(f,x)$
and $(h,y)$ in the graph $\lampp$ if $x,y$ are adjacent in $G$ and
$f(z)=h(z)$ for $z\notin \{x,y\}$.  These wreath products are 
called {\em lamplighter graphs\/} because of the following
interpretation: place a lamp at each vertex of $G$; then a vertex of
$\lampp$ consists of a configuration $f$ indicating which lamps are
on, and a lamplighter located at a vertex $x \in V_G$.  

\begin{figure}
\centering
\epsfig{file=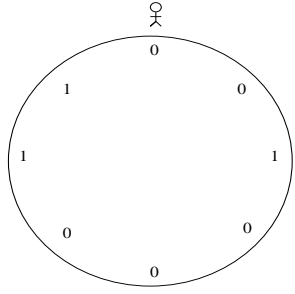, width=1.5 in, height=1.5 in}
\caption{Lamplighter group over a cycle}
\end{figure}
\smallskip

In this paper we estimate mixing time parameters for random walk on
a lamplighter graph $\lampp$ by relating them to hitting and covering
times in $G$.  When $G$
is the two dimensional discrete torus $\Z_n^2$, we prove 
in Theorem \ref{Zdtheorem} 
that: 
\begin{itemize}
\item the relaxation time $\Trel(\lampztwo)$ 
is of order ${n^2 \log n}$ ;
\item
the mixing time in total variation, $\Tve{\lampztwo}$ 
is asymptotic to $ c n^2 \log^2 n$ ; 
\item
the uniform mixing time $\tauve{\lampztwo}$
is asymptotic to $Cn^{4}$.
\end{itemize}
(We recall the  definitions of these mixing time parameters 
 in (\ref{reltimedef})-(\ref{uniformmixdef}).)
The general correspondence between notions of mixing on $\lampp$ and
properties of random walk on $G$ is indicated in the following table:

\begin{center}
{
\def\arraystretch{1.5}
\begin{tabular}{|c|c|} \hline
underlying graph $G$ & lamplighter graph $\lampp$ \\ \hline
maximal hitting time $t^*$ & relaxation time $\Trel$ \\ \hline
expected cover time $\E\C$ & total variation mixing time $\T$ \\
\hline
$\inf\{t: \E2^{|S_t|} < 1+\epsilon \}$ & uniform mixing time
$\tau$ \\ \hline
\end{tabular}
}
\end{center}
where $S_t$ is the set of unvisited sites in $G$ at time $t$.  The
connections indicated in this table are made precise in Theorems
\ref{reltimethm}-\ref{unifconvthm} below.

When $G$ is a Cayley graph, the wreath product $\lampp=\Z_2 \wr G$ 
is the semi-direct product $G \ltimes \sum_G \Z_2$, where 
the action on $\sum_G \Z_2$ is by coordinate shift.  This means that 
multiplication in $\lampp$ is given by $(f,x)(h,y)=(\psi,xy)$, 
where $\psi(i)=f(i)+h(x^{-1} i)$.  
In the case when $G$ is a cycle or a complete graph, random walks 
on $\lampp$ were analyzed by H{\"a}ggstr{\"o}m and Jonasson~\cite{HagJon}.

\medskip
\noindent{\bf Definitions}.
Let $\{X_t\}$ be an irreducible Markov chain on a finite
graph $G$ with transition probabilities given by
$p(x,y)$.  Let $p^t(x,y)$ denote the
$t$-fold transition probabilities and $\mu$ the stationary distribution.
The {\em relaxation time} is given by
\begin{equation}
\label{reltimedef}
\Trel=\max_{i: |\lambda_i|<1} \frac{1}{1-|\lambda_i|}
\end{equation}
where the $\lambda_i$ are the eigenvalues of the transition matrix
$p(x,y)$.
  The
{\em
$\epsilon$-mixing
time in total variation} $\Tve{G}$ and the {\em $\epsilon$-uniform} 
mixing time
$\tauve{G}$ are defined
by:
\begin{equation}
\label{totalvarmixdef}
\Tve{G}=\min \left\{t: \frac{1}{2}
\sum_y |p^t(x,y) - \mu(y) | \leq \epsilon \ \forall \ x \in G \right\}
\end{equation}
and
\begin{equation}
\label{uniformmixdef}
\tauve{G} = \min\left\{t: \left| \frac {p^t(x,y)-\mu(y)}{\mu(y)}
\right| \leq  \epsilon \ \forall \ x, y \in G \right\}.
\end{equation}
When the graph $G$ is clear, we will often abbreviate
\begin{equation}
\T=\Tg{G}=\T(1/(2e),G).
\end{equation}
Another key parameter for us will be the {\em maximal hitting time} 
\begin{equation}
\label{maxhitdef}
t^*=t^*(G)=\max_{x,y} \E_x T_y,
\end{equation}
where $T_y$ is the hitting time of $y$.


The random walk we analyze on $\Z_2 \wr G$ is constructed from a
random walk on $G$ as follows.  Let $p$ denote the transition
probabilities in the wreath product and $q$ the transition
probabilities in $G$.  For $a\neq b$, $p[(f,a),(h,b)]=q(a,b)/4$ if $f$ and 
$h$
agree outside of $\{a,b\}$, and when $a=b$, $p[(f,a), h(a)]=q(a,a)/2$ if 
$f$ and $h$ agree off of $\{a\}$.
A more intuitive description of this is
to say that at each time step, the current lamp is randomized, the
lamplighter moves, and then the new lamp is also randomized.  The
second lamp at $b$ is randomized in order to make the chain
reversible.  To avoid periodicity problems, we will assume that the
underlying random walk on $G$ is already aperiodic.

Our first theorem describes the mixing time of the random walk on the
wreath product when the lamplighter moves in the
$d$-dimensional discrete torus $\Z_n^d$.
\begin{theorem}
\label{Zdtheorem}
For the random walk $\{X_t\}$ on $\lampztwo=\Z_2 \wr \Z_n^2$ in which
the lamplighter performs simple random walk with holding probability
$1/2$ on $\Z_n^2$, the relaxation time satisfies
\begin{equation}
\label{Ztworeltime}
\frac{1}{\pi \log 2} \leq \frac{\Trel(\lampztwo)}{n^2 \log n} \leq
\frac{16}{ \pi \log 2} + o(1).
\end{equation}
For any $\epsilon>0$, the total variation mixing time satisfies
\begin{equation}
\label{Ztwotv}
\lim_{n\rightarrow \infty}
\frac{\Tve{\lampztwo} }{n^2 \log^2 n} =\frac{8}{\pi},
\end{equation}
and the uniform mixing time satisfies
\begin{equation}
\label{Ztwounif}
C_2 \leq \frac{\tauve{\lampztwo}}{n^4} \leq C_2^\prime
\end{equation}
for some constants $C_2$ and $C_2^\prime$.  The uniform mixing
time also has a sharp threshold and
\begin{equation}
\label{Ztwounifsharp}
\taug{\lampztwo}-\tauve{\lampztwo} = O(n^2 \log n) \, .
\end{equation}
More generally,
for any dimension $d\geq 3$, there are constants $C_d$
and $R_d$ independent
of $\epsilon$ such that
on $\Z_2 \wr \Z_n^d=\lampzd$,
the relaxation time satisfies
\begin{equation}
\label{Zdreltime}
\frac{R_d}{4 \log 2} \leq \frac{\Trel(\lampzd)}{n^d} \leq \frac{8
R_d}{\log 2} + o(1),
\end{equation}
the total variation mixing time satisfies
\begin{equation}
\label{Zdtv}
\frac{R_d}{2} + o(1) \leq \frac{\Tve{\lampzd}}{n^d \log n} \leq
R_d +o(1),
\end{equation}
and the uniform mixing time satisfies
\begin{equation}
\label{Zdunif}
C_d \leq \frac{\tauve{\lampzd}}{ n^{d+2}} \leq C_d^\prime \, ,
\end{equation}
with a sharp threshold in that
\begin{equation}
\label{Zdunifsharp}
\taug{\lampzd}-\tauve{\lampzd} =O(n^d) \,.
\end{equation}
\end{theorem}

The parameter $R_d$ is the expected number of returns to 0
by a simple random walk in $\Z^d$, and is
given by equation (75) in Chapter 5 of
\cite{AldFil}.
The reason why the exact limit in (\ref{Ztwotv}) can be computed is
related to the fractal structure of the unvisited set at times up to the
covering time.  The geometry of this set is not sufficiently well
understood in higher dimensions to make it possible to eliminate the
factor of two difference between the upper and lower bounds of
(\ref{Zdtv}).

\noindent
{\bf Remark.}
  In one dimension,
the total variation mixing time on $\Z_2 \wr \Z_n$
was studied in
\cite{HagJon}, and shown to be on the order of $n^2$.
More generally,
the case of $G \wr \Z_n$ for any finite $G$ has similar behavior
\cite{Reyes}.
For the walks we consider, it is easy to show that the relaxation time in
one dimension is on the order of $n^2$, while the uniform mixing time is
on the order of $n^3$.  We compute the asymptotic constant for the 
uniform
mixing time in Section \ref{theoremone}.

The $d$ dimensional result for uniform mixing is
suggested by analogy with the infinite case of $\Z_2 \wr \Z^d$,
where it takes $n^{d+2}$ steps for the probability of
being at the identity to become $\exp(-n^d)$, see \cite{PitSC99b}.

The proofs of Theorem \ref{Zdtheorem} differ in dimensions 1 and 2 from
higher dimensions, partly because a random walk on $\Z^d$ is transient in
dimensions 3 and above.  Many of the ideas for the proof of the higher
dimensional case are actually much more general than what is necessary
for the torus.

\begin{theorem}
\label{reltimethm}
Let $\{G_n\}$ be a sequence of vertex transitive graphs and let $\lamp$ 
denote
$Z_2 \wr G_n$.  As $|G_n|$ goes to infinity,
\begin{equation}
\label{reltimeeq}
\frac{1}{8 \log 2} \leq \frac{\Trel(\lamp)}{t^*(G_n)} \leq
\frac{2}{\log 2}+ o(1).
\end{equation}
\end{theorem}

The lower bound in (\ref{reltimeeq}) is proved using the variational
formula for relaxation time, and the upper 
bound uses a coupling argument that was introduced in \cite{Chen98} 
(see also \cite{BKMP}).  The geometry of lamplighter graphs allows us to 
refine this coupling argument and restrict attention to pairs of
states such that the position of the lamplighter is the same
in both states.

\begin{theorem}
\label{tvconvthm}
Let $\{G_n\}$ be a sequence of vertex transitive graphs with $|G_n|
\rightarrow \infty$, and $\Cn$ denote the cover time for simple random
walk on $G_n$.
For any $\epsilon>0$, there exist constants $c_1, c_2$ depending
on $\epsilon$ such that the total variation mixing time satisfies
\begin{equation}
\label{tvconveq}
[c_1+o(1)] \E\Cn \leq \Tve{\lamp} \leq
[c_2+o(1)] \E\Cn.
\end{equation}
Moreover, if the maximal hitting time satisfies
$t^* =o(\E\Cn)$, then for all $\epsilon>0$,
\begin{equation}
\label{tvconveqsharp}
\left[\frac{1}{2} + o(1)\right] \E\Cn \leq \Tve{\lamp}
\leq [1+o(1)] \E\Cn.
\end{equation}
\end{theorem}

The condition $t^*=o(\E\Cn)$ implies that the cover time has a sharp
threshold, that is $\Cn/\E\Cn$ tends to 1 in probability \cite{Ald}.
Theorem \ref{tvconvthm} thus says that in situations that give a sharp
threshold for the cover time of $G_n$, there is also a threshold for the
total variation mixing time on $\lamp$, although the factor of 2
difference between the bounds means that we have not proved a sharp
threshold.  When $\Teg{G_n}=o(|G_n|)$, (\ref{tvconveq})
is Corollary 2.13 of \cite{Reyes}.

The different upper and lower bounds in
(\ref{tvconveqsharp})
cannot be improved without further hypotheses,
as the limit exists and is equal to the lower bound
when $G_n$ is the complete graph $K_n$ (see \cite{HagJon}) and
the upper bound when $G_n=\Z_n^2$ (Theorem \ref{Zdtheorem}).
The reason for this difference has to do with the geometry of the 
last points that are visited, which are uniformly distributed on $K_n$ but 
very far
from uniformly distributed on $\Z_n^2$.  This difference will be 
emphasized again later on.
We will see
later that the
condition that the $G_n$ are vertex transitive can be replaced
by a condition on the cover times (Theorem \ref{Matthewsthm})
or on $t^*$ (Theorem \ref{tvonrapidmixthm}).

\begin{theorem}
\label{unifconvthm}
Let $\{G_n\}$ be a sequence of regular graphs for which
$|G_n|\rightarrow \infty$ and
the maximal hitting time satisfies
$t^* \leq K |G_n|$ for
some constant $K$.  Then there are constants
$c_1,c_2$ depending on $\epsilon$ and $K$ such that
\begin{equation}
\label{unifconvbnd}
c_1 |G_n|( \Trel(G_n)+ \log|G_n|)  \leq
\tauve{\lamp} \leq c_2 |G_n| (\Tg{G_n}+ \log |G_n| ).
\end{equation}
\end{theorem}

\begin{theorem}
\label{sharpthm}
Let $\{G_n\}$ be a sequence of vertex transitive graph such that
$|G_n|\rightarrow \infty$.  Then the uniform mixing time
$\tauve{\lamp}$ has a sharp
threshold in the sense that for all $\epsilon>0$,
\begin{equation}
\label{sharpeqn}
\lim_{n\rightarrow \infty}\frac{\tauve{\lamp}}{\taug{\lamp}}=1.
\end{equation}
\end{theorem}


The intuition behind both Theorems \ref{tvconvthm} and \ref{unifconvthm}
can be most easily explained by considering the case when $G_n$ is the
complete graph $K_n$, with  a loop added at each
vertex of $K_n$.
The position $\pi(X_t)$ of the lamplighter then performs a simple
random walk on $K_n$ with holding probability $1/n$, and every
lamp that is visited is randomized. 
Thus the (random) cover time for the
walk on $K_n$, which is sharply concentrated near $n\log n$,
is a strong uniform time (see \cite{AldDia})
for the configuration of the lamps.
A strong uniform time (which bounds the mixing time in total
variation)
for the walk on the wreath product is obtained
by adding one more step to randomize the location of the lamplighter.

A special property of the complete graph $K_n$ is that the unvisited set 
is
uniformly distributed.  As shown in \cite{HagJon}, what is actually
needed for mixing in total variation is for the size of the uncovered set 
to be $O(\sqrt{n})$.
At that time, the central limit theorem fluctuations in the number
of lamps that are on are of at least the same order
as the number of unvisited sites. The amount of time needed for this to 
happen is $(n \log n)/2$ steps, resulting in
a sharp phase transition after $(n\log n)/2$ steps.  The factor of
$2$ difference between
the upper and lower bounds in (\ref{tvconveqsharp}) comes from the 
question of
whether or not it suffices to cover all but the last $\sqrt{n}$ sites of
the graph.
For many graphs, the amount of time to cover all but the last 
$\sqrt{|G_n|}$ sites is $\E\Cn/2$, which will be the lower bound of 
(\ref{tvconveqsharp}).  
When the unvisited sites are clumped together instead of being
uniformly distributed, it will turn out to be necessary to visit all the
sites, and the upper bound of (\ref{tvconveqsharp}) will be sharp.

Convergence in the uniform metric depends upon the moment generating
function
of the cover time rather than the mean.
Let $R_t$ denote the set of lamps that the lamplighter has visited by
time $t$, and $S_t=K_n\setminus R_t$ the set of unvisited lamps.
Because any lamp that has been visited is randomized,
the probability that the lamps are in
any given configuration at time $t$ is at most $\E2^{-|R_t|}$.
To
get convergence in the uniform metric, we need $\E2^{-|R_t|}$ to
be within a factor of $(1+\epsilon)$ of $2^{-n}$, or, equivalently, to
have $\E2^{|S_t|}\leq 1+ \epsilon.$

        Considering the Markov
process in continuous time, the visits to any given site in
$K_n$ are independent Poisson processes with rate $n^{-1}$,
so the probability that a given site has not been visited
by time $t$ is $\exp(-t/n)$. Therefore $\E2^{|S_t|}=\E \prod 2^{I_i(t)}$,
where $I_i(t)$ are indicator functions for the event that
 site $i$ has not been visited.  Since
$\E2^{I_i(t)}=[1+\exp(-t/n)]$, we have
$$
\E2^{|S_t|}= [1+\exp(-t/n)]^{n-1}.
$$
For $t=n \log n + cn$, this gives
$$
\E2^{|S_t|} = \exp\left[o(1)+ e^{-c} \right],
$$
so there is a sharp threshold at time $n\log n$.

Since the position of the lamplighter is randomized after
each step, this indicates that $\tauve{K_n^\diamondsuit}$ has a
sharp threshold at $n \log n$ as well.
In this example, the upper and lower bounds of
(\ref{unifconvbnd}) are of the  same order of magnitude and thus are
accurate up to constants.

\section{More examples}
\label{examples}

{ \bf Example (Hypercube)}
On the hypercube $\Z_2^n$, the maximal hitting time is on the order of
$2^n$.
The cover time is on the order of $n2^n$, but
$\Trel(\Z_2^n) =n$ and $\Tve{\Z_2^n} \sim (n \log n)/2$
(\cite{AldFil} Chapter 5 Example 15).
By Theorem \ref{reltimethm}, $\Trel(\Z_2 \wr \Z_2^n)$ is on
the order of $2^n$, and
Theorem \ref{tvconvthm} shows that the convergence time in total
variation on $\Z_2\wr \Z_2^n$ is on the order of $n2^n$.  For uniform
convergence,
the upper and lower bounds of Theorem \ref{unifconvthm} differ, giving,
up to constants, $n2^n$ and $(n \log n) 2^n$.
While Theorem \ref{sharpthm} shows that there is a sharp threshold,
the location of that threshold is unknown for this example.  Thus, while
the relaxation time is less than the mixing times, we do not know whether
or not the total variation and uniform mixing times are comparable.

{\bf Example (Expander graphs)} For a graph $G_n$, the conductance
$\Phi_n$ is the minimum ratio of edge boundary to volume of all
sets $S\subset G_n$ such that $|S|\leq |G_n|/2$.  A sequence of graphs
$\{G_n\}$ is called a sequence of {\em expander graphs} if there
exists a $\delta>0$ such that $\Phi_n \geq \delta$ for all $n$.

For expander graphs, the maximal hitting time on $G_n$ is on the order
of $|G_n|$, while the cover time of $G_n$ is of order $|G_n| \log
|G_n|$ because $\Trel$ is bounded (see \cite{BroKar}).  Theorem
\ref{tvonrapidmixthm} shows that $\Tve{\lamp}$ is of order $|G_n|
\log |G_n|$.  The mixing time on such graphs is of  order $\log
|G_n|$, and so Theorem \ref{unifconvthm} shows that $\tauve{\lamp}
\leq C(\epsilon) |G_n| \log |G_n|$.  Since $\tauve{G} \geq \Tve{G}$
for any $G$, the matching lower bound holds as well.  As a result,
$\{G_n\}$ is a sequence of graphs of bounded degree such that the
total variation and uniform mixing times for the lamplighter graphs
are of the same order.

\section{Relaxation time bounds}

In this section we prove Theorem \ref{reltimethm}.  For the lower bound,
we will use the variational formula for the second eigenvalue.
Let $|\lambda_2|$ denote the magnitude of the second largest eigenvalue 
in
absolute value of a transition matrix $p$, and let $\pi$ be the stationary 
distribution of $p$.
The variational characterization for the eigenvalue says that
\begin{equation}
\label{variationalformeq}
1-|\lambda_2|=\min_{\varphi: \var \varphi >0} 
\frac{\dir{\varphi}{\varphi}}{\var 
\varphi},
\end{equation}
where the Dirichlet form $\dir{\varphi}{\varphi}$ is given by
\begin{align}
\label{dirichletformeq}
\dir{\varphi}{\varphi} & =\frac{1}{2} \sum_{x,y} 
[\varphi(x)-\varphi(y)]^2 \pi(x) p(x,y)\\
& = \frac 12 \E \left[\varphi(\xi_1)-\varphi(\xi_0)\right]^2
\end{align}
and $\{\xi_t\}$ is a Markov chain started in the stationary 
distribution. 
For the lower bound of (\ref{reltimeeq}),
we
use (\ref{variationalformeq}) to show that the
spectral gap for the transition kernel
$p^{t}$ is bounded away from $1$ when $t=t^*/4$.  Fix a vertex $o\in G$,
and let $\varphi(f,x)=f(o)$. Then $\var \varphi=1/4$
   and by running for $t$ steps,
$$
\dirt{\varphi}{\varphi}
= \frac 12 \E \left[\varphi(\xi_t) -\varphi(\xi_0)\right]^2
= \frac 12 \sum_{x\in G} \nu(x) \frac 12 \P_x[T_o <t] ,
$$
where 
$\nu$ is the
stationary measure on $G$, and $\{\xi_t\}$ is the stationary Markov
chain on $\lampp$.  For 
any $t$,
$$
\E_x T_o \leq t + t^* (1-\P_x[T_o<t]) .
$$
For a vertex
transitive graph, we have by  Lemma 15 in Chapter 3
of \cite{AldFil}, that  
$$t^*\leq 2 \sum_{x\in G} \nu(x) \E_x T_0 .
$$
Let $\E_\nu=\sum_x \nu(x) \E_x$
and $\P_\nu=\sum_x \nu(x) \P(x)$.  Then
$$
t^* \leq 2 \E_\nu T_o \leq 2t + 2t^* [1-\P_\nu(T_o<t)].
$$
Substituting
$t=t^*/4$ yields
$$
\P_\nu [T_0 < t^*/4] \leq \frac 34.
$$
We thus have
$$1-|\lambda_2|^{t^*/4} \leq \frac 34,$$
and so
$$
\log 4 \geq \frac{t^*}{4}(1-|\lambda_2|),
$$
which gives the claimed lower bound on $\Trel(\lampp)$.

For the upper bound, we use a coupling argument from \cite{Chen98}.
Suppose that $\varphi$ is an eigenfunction for $p$
with eigenvalue $\lambda_2$.  To conclude that $\Trel({\lampp})
\leq  \frac{(2+o(1)) t^*}{\log 2}$, it suffices to show that
$\lambda_2^{2t^*}\leq 1/2$.  For a configuration $h$ on $G$,
let $|h|$ denote the Hamming length of $h$.  Let
$$
M=\sup_{f,g,x} \frac{|\varphi(f,x)-\varphi(g,x)|}{|f-g|}
$$
be the maximal amount that $\varphi$ can vary over two elements
with the same lamplighter position.  If $M=0$, then
$\varphi(f,x)$ depends only on $x$, and so $\psi(x)=\varphi(f,x)$
is an eigenfunction for the transition operator on $G$.  Since
$\Trel(G) \leq t^*$ (see \cite{AldFil}, Chapter 4), this would
imply that $|\lambda_2^{2t^*}| \leq e^{-4}$.  We may thus assume that
$M>0$.

Consider two walks, one started at $(f,x)$ and one at
$(g,x)$.  Couple the lamplighter component
of each walk and adjust the configurations
to agree at each site visited by the lamplighter.  Let
$(f^\prime,x^\prime)$ and $(g^\prime, x^\prime)$ denote
the position of the coupled walks after $2t^*$ steps.
Let $K$ denote the transition operator of this coupling.
Because $\varphi$ is an eigenfunction,
\begin{align*}
\lambda_2^{2t^*} M & = \sup_{f,g,x} \frac{ |p^{2t^*} \varphi(f,x)-
p^{2t^*}  \varphi(g,x)|}{|f-g|}\\
& \leq \sup_{f,g,x} \sum_{f^\prime, g^\prime, x^\prime}
K^{2t^*}[(f,g,x)\rightarrow (f^\prime, g^\prime, x^\prime)]
\frac{|\varphi(f^\prime, x^\prime)-\varphi(g^\prime, x^\prime)}
{|f^\prime-g^\prime|}\frac{|f^\prime-g^\prime|}{|f-g|}\\
& \leq M \sup_{f,g,x} \frac{\E|f^\prime-g^\prime|}{|f-g|}.
\end{align*}

But at time $2t^*$, each lamp that contributes to $|f-g|$
has probability of at least $1/2$ of having been visited, and
so $\E|f^\prime-g^\prime|\leq |f-g|/2$.  Dividing by $M$
gives the required
bound of $\lambda_2^{2t^*}\leq 1/2$.

\section{Total variation mixing}

In this section, we
prove Theorem \ref{tvconvthm} and then some generalizations.  The total
variation claims of Theorem \ref{Zdtheorem} except for the sharp constant
in
(\ref{Ztwotv}) are applications of Theorem \ref{tvconvthm}, and we defer 
the
computation of the sharp constant in two dimensions until later.

As was mentioned in the discussion of $\Z_2 \wr K_n$, the key to the lower 
bounds comes from running the random walk until the lamplighter visits all 
but $\sqrt{|G_n|}$ lamps, which takes time $\E\Cn/2$ when $t^*=o(\E\Cn)$.

\begin{proof}[Proof of Theorem \ref{tvconvthm}]
We will first prove the lower bound
of (\ref{tvconveqsharp}) since the proof is simpler
when $t^*=o(\E\Cn)$.
Let $S\subset \Z_2 \wr G_n$ be the set
of elements $(f,x)$ such that the configuration $f$ contains
more than $|G_n|/2 + K |G_n|^{1/2}$ zeroes.  Fix $\epsilon \in (0,1)$.
For sufficiently large
$K$ and $n$, we have
      $\mu(S)\leq (1-\epsilon)/4$.  
Fix a basepoint $o\in G_n$, and let $id=(\0,o)$ denote
the element of $\lampp$ corresponding to all the lamps being off and the 
lamplighter being at $o$.
The claim is that for large $n$,
at time  $t_n=\E\Cn/2$,
\begin{equation}
\label{tvlowerneed}
p^{t_n}(id,S)\geq \frac{1+\epsilon}{4}
\end{equation}
and thus $\Tvef{\lamp}\geq \E\Cn/2$.
To see this, let $\Cn(k)$ be the first time that the lamplighter has 
visited all but $k$ lamps, and let
$m_n$ be given by
$$m_n=\min\{t: \P(\Cn(K|G_n|^{1/2}) >t) \leq (1+\epsilon)/2\}.$$
Let
$0<s_1<s_2<s_3<s_4$
be stopping times such that $s_1$ and $s_3-s_2$ are minimal strong uniform 
times, 
and $s_2-s_1$ and $s_4-s_3$ are times to hit all but $K|G_n|^{1/2}$ sites.  
As a result, at time $s_1$ the walk $X_n$ is uniformly distributed, the 
walk covers all but
$K|G_n|^{1/2}$ points of $G_n$ on both intervals $[s_1,s_2)$ and
$[s_3,s_4)$, and conditioned on $X_{s_2}$, the walk is uniformly
distributed at time $s_3$.  

For any finite graph $G$,  denote $\T=\Teg{G}$.
 let $q(\cdot,\cdot)$ be the transition kernel for 
a random walk on $G$ with  stationary distribution $\nu$.
Then the separation distance 
\begin{equation} \label{sepp} 
1-\min_y\{q^{4 \T}(x,y)/\nu(y)\} 
\end{equation}
at time $4\T$, is at most $1/e$
 (see \cite{AldDia} or \cite{AldFil}, Chapter 4, Lemma 7).

Returning to $G_n$, since $\P[s_1>t]$ is the 
separation distance at time $t$ by \cite{AldDia}, we have 
$$\P[s_1> 4 \Teg{G_n}] \leq 1/e .$$
Thus with probability bounded away from $0$, we 
have $$s_4 \leq 2(m_n +4\Teg{G_n}).
$$
But at time $s_4$, the probability of a given point being uncovered is at 
most
$$(K|G_n|^{-1/2})^2,$$
and so the expected number of uncovered points at
time $s_4$ is at most $K^2$.
In particular, with probability at least $1/2$, there are fewer than
$2K^2$ uncovered sites.  Each additional run of $2t^*$ steps has
probability at least $1/2$ of hitting one of these sites, so
with probability $1/4$, all sites are covered after $8K^2 t^*$ more
steps.
Therefore the probability of
having covered the entire space by time $2(m_n+\Teg{G_n})+8K^2t^*$ is
bounded away from 0.  When $t^*=o(\E\Cn)$,
Aldous \cite{Ald} showed that
$\Cn/\E\Cn \rightarrow 1$
in probability, so
$$
\E\Cn \leq (1+o(1)) \left[ 2\Teg{G_n} +2 m_n + 8K^2 t^* \right].
$$
Since $t^*$ and thus $\Teg{G_n}$ are much smaller than $\E\Cn$, this
means
that
$m_n\geq \E\Cn(1+o(1))/2$.
     But at time $m_n-1$, the probability of having
$K|G_n|^{1/2}$ uncovered points is greater than $(1+\epsilon)/2$,
and so the probability of having at least
$*|G_n|+ K |G_n|^{1/2})/2$ zeroes in the configuration
at time $m_n$ tends to $(1+\epsilon)/4$.  This proves the lower bound of 
(\ref{tvconveqsharp}).

For the lower bound of (\ref{tvconveq}), we iterate the
above process one more time to get better control of the probabilities.
More precisely, let
$$
r_n=\min
\left\{t: \P\left(\Cn(|G_n|^{7/12}) > t\right)\leq (1+\epsilon)/2
\right\}.
$$
Let $\P_w$ denote the probability measure for the random walk starting at 
$w$.
By running for
$4 \Teg{G_n}+ r_n$ steps 3 times, a similar argument as before shows that 
the expected number of unvisited sites is small, meaning that
$$
\P_w\left[\Cn> 3(4 \Teg{G_n} +r_n)\right]\leq c \frac{(1+\epsilon)^3}{8}
\left(1- n^{-5/2} \right)^3.
$$
But if $\P_w(C_n>t)\leq x$ for all $w$ then $\E \Cn \leq t(1-x)^{-1}$ 
simply
by running the chain for intervals of length $t$ until all
of $G_n$ is covered in one interval.
We therefore have
$$
\inf\left\{t: \P(\Cn>t) \leq c \frac{(1+\epsilon)^3}{8}\right\}
\geq c_1(\epsilon) \E\Cn
$$
for some constant $c_1(\epsilon)$.  Since $t^*$ (and thus $\Teg{G_n}$ by
\cite{AldFil}, Chapter 4 Theorem 6) is $o(\E\Cn)$,
this shows that $r_n\geq c_1(\epsilon) \E\Cn
(1+o(1))$.  For large $n$, however, the total variation distance on
$\lamp$
is still at least $\epsilon$ at time $r_n$, which proves
the lower bound of (\ref{tvconveq}).

For the upper bound, let $\pi$ denote the projection from $\lampp$ to $G$
given by the position of the lamplighter.
Let $\mu$ denote the stationary measure on $\lampp$ and $\nu$ the
stationary measure on $G$.  Note that $\mu(f,x)=2^{-|G|}\nu(x)$.
By the strong Markov property, at time
$t=\C+1+k$ we have $p^t[id,(f,x)]=2^{-|G|}
\P_{X_\C}[\pi(X_{k+1})=x]$.  For
any $k>0$,
\begin{multline}
\label{tvdecomp}
\sum_{(f,x)\in \Z_2 \wr G} \bigg| p^t[id,(f,x)]-\mu(f,x)\bigg|
\leq \P(t \leq \C + k ) \\
+ \sum_{x\in G} \bigg| \P_{X_\C}[\pi(X_{k+1})=x]-\nu(x) \bigg|.
\end{multline}
Letting $t_1(\epsilon)=\min\{t: \P(\C>t) \leq \epsilon/2\}$,
we see that
$\Tve{\lamp}$
is at most $t_1(\epsilon)+\Tvef{G_n}$.  The desired upper
bound comes from the fact that $t_1 \leq [c_2(\epsilon)+o(1)] \E\Cn$
and $\Tvef{G_n}$
is of a lower order than $\E\Cn$.  When $t^*=o(\E\Cn)$, the fact that
$\Cn/\E\Cn$
tends to $1$ in probability implies that $c_2(\epsilon)=1$ works
for all $\epsilon$.
     \end{proof}

In many situations, the assumption of vertex transitivity
in Theorem \ref{tvconvthm} can be relaxed.

\begin{theorem}
\label{Matthewsthm}
Suppose that $G_n$ is a sequence of graphs for which simple random walk
satisfies
$$\lim_{n\rightarrow \infty} \frac{t^* \log|G_n| }{ \E\Cn }= 1,$$
and
$|G_n|\rightarrow \infty$.
Then for any $\epsilon>0$,
\begin{equation}
\left[\frac{1}{2}+o(1)\right] \E\Cn \leq \Tve{\lamp} \leq [1+o(1)]
\E\Cn \, .
\end{equation}
\end{theorem}

\begin{proof}
Since the proof of the upper bound for Theorem \ref{tvconvthm} did not
rely on vertex transitivity, we only need to prove the lower bound.  To
do
so, it suffices to show that for any $\delta>0$,
     \begin{equation}
\label{expcovereqn}
\lim_{n\rightarrow \infty} \P[\Cn(|G_n|^\alpha)<(1-\alpha-\delta) 
\E\Cn]= 0 \,,
\end{equation}
and then follow the proof of Theorem \ref{tvconvthm}.

The way that we will prove this is by contradiction.
If $$\limsup_{n\rightarrow \infty}
 \P[\Cn(|G_n|^\alpha)<(1-\alpha-\delta) \E\Cn]>0,$$
 we will show that we can find an $\epsilon>0$ such that
 $$\limsup_{n\rightarrow \infty}
 \P[\Cn<(1-\epsilon) \E\Cn]>0.$$
 But $\Cn/\E\Cn\rightarrow 1$ in probability, so this will be
 a contradiction.  To do this,
let $S_n$ denote the final $|G_n|^\alpha$ points in $G_n$ that are hit by
$\{X_t\}$.  We will first follow Matthews' proof \cite{Mat} to show that
$$\E[\Cn-\Cn(|G_n|^\alpha) \mid S_n] \leq \alpha \E\Cn.$$

Label the points of $S_n$ from $1$ to $|G_n|^\alpha$ with a labelling
chosen uniformly from the set of all such labellings, and let $\tilde
\Cn(k)$ denote the first time after $\Cn(|G_n|^\alpha)$ that the first $k$
elements of $S_n$ (according to our random labelling) are covered by
$\{X_t\}$.  Note that $\P[\tilde \Cn(k+1) > \tilde \Cn(k)]=1/(k+1)$ by
symmetry, and thus $\E[\tilde \Cn(k+1) - \tilde \Cn(k)] \leq
t^*/(k+1)$.
Summing this telescoping series, we see that $\E\Cn-\E\Cn(|G_n|^\alpha)
\leq
t^* [\log (|G|^\alpha)+1] \sim \alpha \E\Cn$.

To complete the proof, note that
\begin{multline*}
\P[\Cn < (1-\epsilon) \E\Cn] \geq \P[\Cn(|G_n|^\alpha)
<(1-\alpha-\delta) \E\Cn]
\\ \times
 \P[\Cn-\Cn(|G_n|^\alpha) < (\alpha 
+\delta -\epsilon) \E\Cn \mid S_n ].
\end{multline*}
For $\epsilon$ such that $\alpha(1+\epsilon)^2 < \alpha +\delta
-\epsilon$, Markov's inequality shows that the final
probability is at least $\epsilon/(1+\epsilon)$ for large
$n$, which completes the proof.
\end{proof}

\begin{theorem}
\label{tvonrapidmixthm}
Let $G_n$ be a sequence of regular graphs such that $\Trel(G_n)\leq
|G_n|^{1-\delta}$ for some $\delta>0$.
Then there are constants $c_1(\epsilon)$ and $c_2(\epsilon)$ such
that
\begin{equation}
\label{tvonrapidmixeq}
\left[c_1+o(1)\right] |G_n|\log |G_n| \leq \Tve{\lamp} \leq [c_2+o(1)]
|G_n| \log |G_n|.
\end{equation}
\end{theorem}

\begin{proof}
This is actually a corollary of a result of Broder and Karlin
\cite{BroKar}.
They showed that $\E\Cn$ is of the order of $|G_n|\log |G_n|$.  Examining
the proof of their lower bound (Theorem 13) shows that they actually
proved the stronger statement that the time to cover all but
$|G_n|^{\alpha+1/2}$ sites is of order $|G_n|\log |G_n|$ for small
enough $\alpha$, from which our result follows as before.
\end{proof}

\section{Uniform mixing}
\label{convproofs}

As discussed in the example of a lamplighter graph
on the complete graph $K_n$, the key to the mixing time on the wreath
product is the number of sites on the underlying
graph that are left uncovered by the projection
of the lamplighter's position.
We thus
first develop some necessary facts about the distribution
of the number of uncovered sites, and
then prove our main results.

Recall from (\ref{sepp}) that
the transition kernel $q$ on $G$ and the stationary measure $\nu$ 
satisfy $q^{4\T}(x,y)\geq (1-e^{-1}) \nu(y)$,
where $\T=\Teg{G}$, 

\begin{lemma}
\label{largedevneed}
Let $\{X_t\}$ be an irreducible Markov chain on a state space 
$G$
and let $\Tret$ denote the first return time to $x$.
Suppose that there are $\epsilon, \delta>0$ such that
$$\P_x(\Tret > \epsilon |G|) \geq \delta>0$$
    for all $x\in G$.
Let $S\subset G$ be a set of $k$ elements.  Then for $k\geq
\T$,
the probability of hitting at least $\delta \epsilon \T/ 4$ elements
of $S$
by time $$4\T+ \epsilon |G| \T k^{-1}$$
is bounded away from 0.
\end{lemma}

\begin{proof}
Let $r=\epsilon G \T k^{-1}$.  For $1 \leq i \leq r$, 
let $I_i$ be an indicator random variable
for the event $\{X_{4\T+i} \in S \}$,
and let $J_i$ be an indicator for the event
that $X_{4\T+i}$ is in $S$ but that the walk does not
return to $X_{4\T+i}$ by time $r$. 
Note that  $I_i\geq J_i$,
and $\sum_1^r J_i$ is the number of distinct elements of $S$
visited between time $4\T+1$ and time $4\T+r$, inclusive.
By running for an initial amount of time $4\T$, the probability of the 
lamplighter being at any given lamp is at least $(1-e^{-1})|G|^{-1} \geq 
(2|G|)^{-1}$.
As a result, $\E J_i \geq \delta k(2|G|)^{-1}$ and thus
$\E\sum_1^r J_i \geq \epsilon \delta \T/ 2$.  To conclude that
$\P(\sum J_i > \epsilon \delta \T/ 4)$ is bounded away from 0, it
suffices
to show that there is a constant $C$ such that $\E(\sum J_i)^2 \leq C
\T^2$.  But
$$\Cov(I_i, I_j) \leq \|I_i\|_2 \|I_j\|_2 \exp[-\Trel^{-1}|i-j|],
$$
and summing the geometric series,
$$\sum_{j=1}^r \Cov(I_i,I_j) \leq \frac{k \Trel}{2 
|G|(1-\exp[-\Trel^{-1}])}.$$ 
Therefore,
$$\sum_{i=1}^r \sum_{j=1}^r \Cov(I_i, I_j) \leq \frac{\delta \epsilon 
\T}{2(1-\exp[-\Trel^{-1}])}.$$
Because $1/2\leq \Trel\leq \T$,
there is a constant $C$ such that
$$\E\Bigl(\sum_{i=1}^r I_i\Bigr)^2 \leq C \Tg{G}^2 \, . 
$$
Since $\E(\sum_{i=1}^r J_i)^2 \leq \E(\sum_{i=1}^r I_i)^2$, this 
completes the proof.
\end{proof}

Lemma \ref{largedevneed} does not apply once we have reduced the number of 
unvisited sites to less than $\T$, but that situation is easier to 
control.

\begin{lemma}
\label{largedevneed2}
Let $\{X_t\}$ be a Markov chain on $G$ such that $t^* \leq K |G|$, and let 
$S \subset G$ be a subset of $s$ elements.  Then the probability of 
hitting at least $s/2$ elements of $S$ by time $2K |G|$ is at least $1/2$.
\end{lemma}

\begin{proof}
Suppose $x\in S$, and let $T_x$ be the first hitting time of $x$.  Since 
$\E T_x\leq t^*\leq K|G|$, we have $\P[T_x \leq 2K|G|] \geq 1/2$.  The 
expected number of elements of $S$ that are hit by time $2K|G|$ is thus at 
least $s/2$, and the number that are hit is at most $s$, so the 
probability of hitting at least $s/2$ is bounded below by $1/2$ as 
claimed.
 \end{proof}

As the next lemma shows,
the hypotheses of Lemma \ref{largedevneed} are easily fulfilled.

\begin{lemma}
\label{transience}
Suppose that $G$ is a regular graph with maximal hitting time
$t^*$.  The distribution of $T_x^+$, the first return time to $x$, 
satisfies
$$\min_{x\in G} \P_x\left(\Tret >\frac{|G|}{2}\right) \geq
\frac{|G|}{2t^*}
       >0.$$
\end{lemma}

\begin{proof}
For random walk on a regular graph,
$\E_x \Tret = |G|$.  But $\E_x \Tret \leq |G|/2 + \P_x(\Tret> |G|/2)
t^*$,
and rearranging terms proves the claim.
\end{proof}

\begin{lemma}
\label{largedevlemma}
Let $\{X_t\}$ be a Markov chain on $G$ such that
$t^* \leq K |G|$.  Let
$S_t$ denote the set of points not covered by time $t$.
Then there are constants $c_1, c_2$, and $c_3>0$ such that if
$t=(1+a)c_1 |G| (\T+ \log|G|)$, then
$\E2^{|S_t|} \leq 1+ c_2 \exp[-a c_3 (\T+\log|G|)]$.
\end{lemma}

\begin{proof}
Let $r=\lfloor |G|\T^{-1} \rfloor$.
For $0\leq i \leq r-1$,
let $k_i=|G|-i \T$, and for $i=r$ let $k_i=0$.
Let $t_i=\min\{t: |S_t| \leq k_i\}$.  For $i<r$, 
Lemmas \ref{largedevneed} and \ref{transience} imply that 
there is a
constant $c=c(K)$ such that
$t_{i}-t_{i-1}$ is stochastically dominated by
$x_i=(c|G| \T k_i^{-1}) Z_i$, where $Z_i$ is a
geometrically distributed
random variable with mean 2.
We will begin by bounding the probability that $t_i$ is much larger than
$|G| (\T+ \log |G|)$ when $i<r$.
By expressing $t_i$ as a telescoping series and using Markov's 
inequality, we see that
for any $\theta \geq 0$,
\begin{align*}
\P[t_i>t] \leq & \P \left[\sum_{j=1}^i x_j > t \right]\\
= & \P\left[ \sum_{j=1}^i \frac{c|G|\T}{k_j} Z_j > t \right]\\
\leq & \exp[-t \theta] \E \exp\left[\sum_{j=1}^i \frac{c |G|\T 
\theta}{k_j} Z_j \right] \\
= & \exp[-t\theta] \prod_{j=1}^i \E \exp \left[ \frac{c|G|T_v \theta}{k_j} 
Z_j \right].
\end{align*}
But for $\alpha \leq 1/3 \leq \log(3/2)$,
$$
\E \exp[\alpha Z_j]= \frac{\exp(\alpha)}{2-\exp(\alpha)} \leq 
\exp(3\alpha)\, .
$$
Let $\theta= k_i(3c|G|\T)^{-1}$.
For $j\leq i$ we have $k_j <k_i$, whence
$$
\frac{c|G|\T \theta}{k_j} \leq \frac 13 \,.
$$
Consequently,
$$\E \exp\left[\frac{c|G|\T \theta}{k_j} Z_j\right] \leq
\exp\left[ \frac{k_i}{k_j}\right].$$
Since $\P[|S_t|>k_i]=P[t_i>t]$, and $i\leq |G| T_v^{-1} +1$, we have
\begin{equation}
\label{upperboundeqn}
\P[|S_t|>k_i] \leq
\exp\left[- t\left(\frac{k_i}{3c\T |G|}\right) + \frac{k_i}{ \T} \log
|G|\right].
\end{equation}
For $i=r$, Lemma \ref{largedevneed2} shows that $t_r-t_{r-1}$ is 
stochastically dominated by a sum of at most $\log_2 (2 \T)$ geometric 
random variables with mean $2 K |G|$.  As a result, 
\begin{equation}
\label{upperboundeqn2}
\P[t_r-t_{r-1}>t] \leq \log_2 (2 \T) \exp\left[-\frac{t}{2K|G|}\right]
\,.
\end{equation}
Breaking the possible values of $|S_t|$ into intervals of length $\T$, we 
get
$$
\E 2^{|S_t|} \leq 1 + \sum_{i=0}^{r} 2^{k_i +\T} \P(|S_t|>k_i).
$$
Let $c_1=\max\{7 c \log 2, 5K \log 2\}$.  For $i<r$, $k_i \geq \T$ and
$$
2^{k_i+\T} \P[|S_t|>k_i] 
 \leq \exp\left[ -\frac 73 (\log_2) k_i a - \frac 13 (\log 2) k_i \right] 
$$
where the second inequality used the fact that $k_i\geq \T$ for $i<r$.
When $i=r$,
\begin{align*}
2^{k_r+\T} \P[|S_t| >k_r] \leq &  2^{\T}[\P(t_{r-1}>t/2) + 
\P(t_r-t_{r-1}>t/2)] \\
 \leq & 2^{T_v} \exp\left[-\frac 73 (1+a)k_{r-1}\log 2\right] \\
& + 2^{\T} \log_2(2 
\T) \exp\left[ -\frac{5 \log 2 (1+a) (\T+ \log |G|)}{4} \right]
\end{align*}
Summing these terms yields the desired bound.
\end{proof}

We will defer the proof of Theorem \ref{Zdtheorem} to Section 
\ref{theoremone},
but will prove here the somewhat more
general Theorem \ref{unifconvthm}.

\begin{proof}[Proof of Theorem \ref{unifconvthm}]
First, we will prove the lower bound on $\tauve{\lamp}$.  The
bound $\tauve{\lamp}\geq c |G_n| \log |G_n|$ is only
relevant when $\Trel(G_n) \leq |G_n|^{1/2}$, in which
case it is a
consequence of Theorem \ref{tvonrapidmixthm} and the general
inequality
 $\tauve{G} \geq \Tve{G}$, applied to $G=\lamp$. 
To establish the lower bound, it remains to prove that
$\tauve{\lamp}\geq c |G_n| \cdot \Trel(G_n)$. 
note that for any $\epsilon>0$, 
we have
$\tauve{\lamp} \geq |G_n|/2$
when $|G_n|$ is large enough.  We may thus assume that 
$\Trel(G_n)$ is bounded away from 1, and for concreteness will assume
$\Trel(G_n)\geq 2$.
Let $\eta$ be a right eigenfunction
of the walk $\{Y_t\}$ on
$G_n$ with eigenvalue $\lambda=1-\Trel^{-1}(G_n)$
and $\sum_g \eta(g)=0$.  Then $\lambda^{-t}
\eta(Y_t)$ is a martingale.  Let $S$ be the smallest subset
consisting of at least half of the vertices of $G_n$ of the form
$S=\{x\in
G_n \, : \, \eta(x) \leq a\}$.  Possibly replacing $\eta$ by $-\eta$, we 
may
assume
that $a\leq 0$.  Let $x_0$ be a vertex of $G$ such that
$\eta(x_0)=M=\max \eta(x)$.

Let $T_S$ denote the hitting time of $S$.
By the optional stopping theorem,
\begin{equation}
\E_{x_0}\left[ \lambda^{-T_S\wedge t} \eta(Y_{T_S \wedge
t})\right]=\E_{x_0} \eta(Y_0).
\end{equation}
Since $a\leq 0$, this yields
\begin{equation}
M \P(T_S > t) \lambda^{-t} \geq M.
\end{equation}
Since $\Trel(G_n)\geq 2$, we have 
$\log \lambda \geq -(2 \log 2) \Trel^{-1}$,
whence $$P(T_S>t) \geq \exp[-(2 \log 2) t \Trel^{-1}].$$  
At time $t$, the
probability that all lamps are off is at least
$P(T_S>t)2^{-|G_n|/2}$, and for
$$t< \Trel(G_n) \cdot \left( \frac {|G_n|}{4} - \frac{\log(1+\epsilon)}
{2 \log 2} \right)$$
the probability that all lamps are off is at
least $(1+\epsilon) 2^{-|G_n|}$, so
$\tauve{\lampp} \geq |G_n| \Trel(G_n) \cdot(1/4+o(1))$.

For the upper bound
$\tauve{\lamp} \leq t= C(\epsilon)|G_n| (\log |G_n|+
\Tve{G_n})$,
we need to prove two things for any $(f,x)\in \lamp$: 

First, that $\P(X_t=(f,x)) \geq 2^{-|G_n|} |G_n|^{-1}(1-\epsilon)$, and
second, that $\P(X_t=(f,x)) \leq 2^{-|G_n|} |G_n|^{-1} (1+\epsilon)$.

For the first bound,
$$\P\bigg(X_t=(f,x)\bigg) \geq |G_n|^{-1} 2^{-|G_n|} \P\bigg(\Cn < t-
\tauvef{G_n} \bigg) (1-\epsilon/2),$$
and the hypotheses of the theorem imply that
$\P(\Cn < t - \tau(\epsilon/2,G_n))$ tends to 1 (see \cite{Ald}).

For the second bound, we first run the walk for
$r=t-\tau(\epsilon/3, G_n)$
steps.
Let $\pi$ denote the projection from $\lamp$ to $G_n$
given by the position of the lamplighter, and
write $X_t=(f_t,\pi(X_t))$.
Since the lamps visited by time $r$ are
all randomized,
$$
\P(f_t=f) \leq \E 2^{|S_r|-|G_n|},$$
 where $S_r$ is the
uncovered set of $G_n$ at time $r$.  
Considering the next $t-r$ steps, we get
$$
\P(\pi(X_t)=x \, | \, X_{j}, j \le r) \leq \max_y \P(\pi(X_{t-r})=y) \,,
$$
 which is at most $(1+\epsilon/3)|G_n|^{-1}$.  

By Lemma \ref{largedevlemma}, $$
\E 2^{|S_r|-|G_n|} \leq (1+\epsilon/2) 2^{-|G_n|}.$$
  Combining these, we have
$$\P\bigg(X_t=(f,x)\bigg)\leq (1+\epsilon) |G_n|^{-1} 2^{-|G_n|} , $$
as required.
\end{proof}

\begin{proof}[Proof of Theorem \ref{sharpthm}]
What we will actually prove is that
\begin{equation}
\label{sharpthmineq}
\tauveg{7\epsilon/8}{\lamp} \leq \tauve{\lamp}
+ 2t^*+ \tauveg{\delta}{G_n}.
\end{equation}
for a suitable $\delta=\delta(\epsilon)$.
Write $\Trel=\Trel(G_n)$.
Given (\ref{sharpthmineq}),
proving the theorem only requires
that $t^*+\tauvef{G_n}=o(\tauve{\lamp})$.
But $$\tauvef{G_n} \leq \log |G_n| \cdot \Trel =o(\tauve{\lamp}),$$
so we only need the condition that $t^*=o(\tauve{\lamp})$.

For graphs $G_n$ such that $t^*/ \Trel\leq |G_n|^{1/2}$, the
lower bound of (\ref{unifconvbnd}) proves the claim.
When $t^*/ \Trel \geq |G_n|^{1/2}$,
Theorem 4 and
Corollary 5 in Chapter 7 of \cite{AldFil}
show that
$$\E\Cn \geq (1/2+o(1)) t^* \log |G_n|.$$
We are thus in a case covered by Theorem \ref{tvconvthm},
and so
$$\tauve{\lamp} \geq \Tve{\lamp} \geq (1/4 + o(1)) t^* \log |G_n|,$$
from which the claim follows.

To show that (\ref{sharpthmineq}) holds,
consider the uncovered set at the times
$r=\tauve{\lamp}$ and $t=r+2t^*$.
Conditioned on the $|S_r|>0$, at time
$t$ at least one more point has been covered with
probability at least $1/2$.  This means that
$$
\E 2^{|S_t|}-1 \leq \frac{3}{4} \left(\E 2^{|S_r|}-1\right).
$$
Because $\E 2^{|S_r|} \leq 1+\epsilon$, this implies that $\E 2^{|S_t|}
\leq
1+
3\epsilon/4$.
Let $u=t+ \tauveg{\delta}{G_n}$.  Then for any $(f,x)\in \lamp$,
\begin{align*}
\P[(f_u, \pi(X_u)) &=(f,x)] \leq \E\left[2^{|S_t|-|G_n|}\right]
\max_{y\in G_n} \P_y [\pi(X_{u-t})=x] \\
& \leq \left( 1 + \frac{3\epsilon}{4}\right) 
(1 + \delta) |G_n|^{-1} 2^{-|G_n|}.
\end{align*}
Taking $\delta$ small enough that
$(1+3\epsilon/4)(1+\delta)<(1+7\epsilon/8)$ proves
that (\ref{sharpthmineq}) holds.
\end{proof}

\section{Proof of Theorem \ref{Zdtheorem}}
\label{theoremone}

We turn now to completing the proof of Theorem \ref{Zdtheorem}.
The statements about relaxation time follow from Theorem
\ref{reltimethm} and standard facts about $t^*$ (see e.g.
\cite{AldFil} Chapter 5).
The mixing results
in dimension $d\geq 3$ follow from Theorems
\ref{tvconvthm} and \ref{unifconvthm}, so only
the dimension $d=2$ results remain to be shown.

\begin{proof}[Proof of Theorem \ref{Zdtheorem}]
As usual, we let $S_t$ denote the set
of lamps that are unvisited at time $t$.

For completeness, we will also show that for $d=1$,
$$
\tauve{\lamp}\sim \frac{64 \log 2}{27 \pi^2}n^3.
$$
For simple random walk with holding probability $1/2$ in 
one dimension, at time $t=\alpha n^3$,
$$
\P(|S_t|=\lambda n) = \exp\left[ -\left(\frac{\pi^2}{4}+o(1)\right)
 \frac{\alpha n}{(1-\lambda)^2} \right]
$$
(see \cite{Spitzer}, Section 21), and so
$$\E 2^{|S_t|} - 1 = n \int_{0}^{1} \exp \left[ -
\left( \frac{\pi^2}{4}+o(1)\right)
\frac{\alpha n}{(1-\lambda)^2} \right] d\lambda.$$
Taking a Taylor expansion about $\rho=1-(\pi^2 \alpha/ 2 \log 2)^{1/3}$ to
estimate the integral, we see that there is a sharp threshold for the
time
at which $\E 2^{|S_t|}=1+\epsilon$ that occurs when
$$
t= \frac{32 \log 2}{27 \pi^2}.
$$
The argument in the proof of Theorem \ref{unifconvthm} shows that this is
also a sharp threshold for $\tauve{\lamp}$.

For $d=2$, we will give a proof that at first seems wasteful but is
surprisingly accurate for the extreme events that are the dominant
contribution to $\E 2^{|S_t|}$.
After $n^2$ steps, there exists an $\alpha>0$
such that a simple random walk on $\Z_n^2$ will cover some ball
of radius $n^\alpha \sqrt 2$ with probability at least
$\delta>0$ \cite{Revesz}.
Divide the
torus $\Z_n^2$ into $n^{2-2\alpha}$ boxes of side length $n^\alpha$.
On time intervals of the form $[2k n^2, (2k+1)n^2]$, there is thus a
probability of at least $\delta$ of the random walk covering at least one
of the $n^{2-2\alpha}$ boxes.  Moreover, because the mixing time on the 
torus is of the order $n^2$,
we see that even conditioned
on which boxes are currently covered, if there are $i$ uncovered boxes
then the walk will cover one of the uncovered boxes with probability at
least $\epsilon i n^{2\alpha-2}$ for some fixed $\epsilon>0$.  
Since the probability of covering a new box in any given interval of 
length $2k^2$ is uniformly bounded away from 
0, the number of such intervals needed to cover a new box is 
stochastically dominated by a geometric random variable.  In particular, 
rescaling implies that
if $t_i=\min\{t: $ there are $i$ covered boxes$ \}$, then for large enough 
$c$
$$
\frac{n^{2-2\alpha}-j}{n^{4-2\alpha}c} (t_{j}-t_{j-1}),
$$ 
is stochastically dominated by a geometric random variable with mean 2.
We will 
use this fact to crudely bound the exponential moment of the number of 
uncovered
sites on the torus.

As before,
\begin{equation}
\label{Ztwounifeq}
\E[2^{|S_t|}] \leq 1+ \sum_{i=1}^{n^{2-2\alpha}} \P(t_i >t) 
2^{n^2-in^{2\alpha}}.
\end{equation}
Arguing as in the proof of Lemma
\ref{largedevlemma},
\begin{equation}
\P(t_i>t) \leq  e^{-t \theta} \E \exp
\left[\sum_{j=1}^{i}
\frac{\theta n^{4-2\alpha} c}{n^{2-2\alpha}-j} Z_j\right]
\end{equation}
and so for $\theta=(n^{2-2\alpha}-i)(3n^{4-2\alpha} c)^{-1}$,
       \begin{equation}
\P(t_i>t)
\leq \exp\left[ -t \frac{n^{2-2\alpha}-i}{3 n^{4-2\alpha} c}
+ (n^{2-2\alpha}-i) \log n \right].
\end{equation}
For $\delta>0$ and $t=(1+\delta) n^4 4c \log 2$, substituting into
(\ref{Ztwounifeq}) shows that the order of $n^4$ steps suffice to 
reduce
$\E 2^{|S_t|}$ to $(1+\epsilon)$
and are thus an upper bound for $\tauve{\lamp}$.

For a lower bound on $\tauve{\lamp}$, let $S$ be the set
of elements in $\lamp$
for which the absolute value of the horizontal component of
the
lamplighter's position is at most $n/4$. The size of $S$ is $n^2/2$, and
the probability of remaining in $S$ for the first $cn^4$ steps
decays like $\exp[-c \alpha n^2]$, and so for small enough $c$, we have
$\E 2^{|S_t|} \geq 2^{n^2/2} \exp[-c \alpha n^2] > 1+\epsilon$.
The sharp threshold in (\ref{Ztwounifsharp}) is due to Theorem 
\ref{sharpthm}.

To show that the sharp threshold claimed in (\ref{Ztwotv})
actually exists, we need to prove that a lower bound for total variation
convergence is $\E\Cn$. Let $\alpha \in (0,1)$.
By \cite{DPRZ}, Sec.\ 4, there exists $\beta=\beta(\alpha)>0$
such that at time $(1-\alpha) \E\Cn$, the uncovered set
contains a ball of radius $n^{\beta}$ with high
probability.
The configuration $f$ at time $(1-\alpha) \E \Cn$ is thus,
with high probability,
identically zero on a ball of radius $(1/3)n^{\alpha/2}$.
But in the stationary measure $\mu$ on $\lamp$, the expected
number of balls of radius $r$ for which the configuration is
identically 0 on the ball is less than $n^2 2^{-r^2}$.  Therefore the
probability of
having a ball of 0's with
radius of a greater order than $(\log n)^{1/2}$ tends
to 0 as $n$ grows.
\end{proof}

\section{Comments and questions}
Random walks on another type of wreath product were analyzed by 
Schoolfield \cite{School} and by Fill and Schoolfield
\cite{FillSchool}. When 
$G$ is the symmetric group $S_n$, instead of letting $S_n$ act on itself 
by left-multiplication, it is perhaps more natural to let $S_n$ act on 
$\Z_n$ by permutation.  As such, $G\wr S_n$ is often used to describe $S_n 
\ltimes \sum_{\Z_n} \Z_2$, while our description yields the Cayley graph 
of $S_n \ltimes \sum_{S_n} \Z_2$.  Mixing times in both total variation 
and $L^2$ norm are carefully analyzed for this alternative description of 
$G\wr S_n$ in \cite{School} and \cite{FillSchool}.

\medskip

Although the results of this paper give a good description
of the mixing on finite lamplighter groups in both
total variation norm and the uniform norm, the
discrepancies in the upper and lower bounds raise
a number of natural questions.

\begin{enumerate}
\item
For the torus in dimension $d\geq 3$, is $\Tve{\lamp}$ asymptotic to
$\E\Cn$ or $\E\Cn/2$?  This question can be stated without reference to
the lamplighter: let $R_t$ be the set of points that simple random walk 
on the torus $G_n$ visits by time t.
Suppose that $\mu$ is uniformly distributed on
$\{0,1\}^{G_n}$ and that $\nu_t$ is uniformly distributed on
$\{0,1\}^{R_t} \times \{0\}^{G_n \setminus R_t}$.  For $1/2 < \alpha <1$, 
what is
$$
\lim_{n\rightarrow \infty} \|\mu -\nu_{\alpha \E \Cn}\|_{TV}?
$$
Although the last $n^{d/2}$ points
to be covered on $\Z_n^d$ are not exactly uniformly distributed,
it is possible that the difference from uniform is dominated by the noise
in the uniform measure on $\{0,1\}^{G_n}$. Proving
that this total variation distance tends to 0 for $\alpha\in(1/2,1)$
    would be a way
    of quantifying that the last points to be visited are nearly uniformly
distributed.
    \item
More generally, for which sequences of graphs $G_n$
is $\Tve{\lamp}$ asymptotic to $\E\Cn$, and
for which is it $\E\Cn/2$?  Can it be asymptotic
to $\alpha \E\Cn$ for some $\alpha \in (1/2,1)$?
\item
Is $\E\Cn$ the right order of magnitude
for the total variation mixing even when the graphs $G_n$ are not vertex
transitive?  That is to say, if $t_n=o(\E\Cn)$, is
$\Tve{\lamp} > t_n$ for large enough $n$?
\item
For $G_n=\Z_2^n$, is the correct order
of magnitude for $\tauve{\lamp}$ equal
to $n 2^n$ or $n (\log n) 2^n$?
\item
Can the upper bound on the uniform mixing time in Theorem 
\ref{unifconvthm} be
replaced by $C |G_n|(\Trel+ \log|G_n|)$?
\item
To what extent are these results generator dependent?  Another type of
walk that is often considered on lamplighter groups is one in which
at each step either
the current lamp is adjusted while the lamplighter holds, or no lamps are
adjusted while the lamplighter moves to a neighbor.  By
comparing Dirichlet forms, this only changes the relaxation time up to
constants, but it is not clear what happens to the mixing times.
Using log-Sobolev inequalities shows that the uniform mixing time can 
change by at most a factor of $\log \log |\lampp| \simeq \log |G|$.  
\end{enumerate}

\noindent{\bf Acknowledgement}
We thank G\'abor Pete for useful comments.

\bibliographystyle{amsplain}
\bibliography{../refs}

\end{document}